\documentclass[12pt]{article}

\usepackage{amssymb}

\usepackage{mathrsfs}

\usepackage{geometry}          

\usepackage{epstopdf}
\usepackage{amsthm}
\usepackage{enumerate}

\usepackage{dsfont}

\newtheorem{theorem}{Theorem}

\newtheorem{proposition}[theorem]{Proposition}

\def\qed{\unskip\nobreak\hfill\penalty50\hskip 3pt\hbox{}\nobreak
\hfill\hbox{\vrule width 4 pt height 10 pt}}

\begin{document}
\bibliographystyle{plain}

\title{A note on the limiting entry and return times distributions for  induced maps}
\author{Nicolai Haydn
\thanks{Mathematics Department, USC,
Los Angeles, 90089-1113. E-mail: $<$nhaydn@math.usc.edu$>$.}}

\maketitle

\abstract{For ergodic measures we consider the return and entry times for a measure
preserving transformation and its induced map on a positive measure subset. 
We then show that the limiting entry and return times distributions are the same for 
the induced maps as for the map on the entire system. The only assumptions 
needed are ergodicity and that the measures of the sets along which the limit is 
taken go to zero.}

\section{Introduction}

In \cite{BSTV} it was shown that on manifolds the limiting return times distributions for an
ergodic map is the same for the induced map on a subset if the limit is along 
a sequence of metric balls. The theorem was proven for measures that satisfy
the Lebesgue density theorem, which include absolutely continuous measures
and also Radon measures. Here we will show a general version of that 
statement that only requires ergodicity and does not impose any other restrictions.

The statistics of entry and return times has been studied to quite some degree in
the last two decades in particular. A number of results have been achieved for 
a variety of systems that have good mixing properties as for instance for 
Axiom A systems or subshifts of finite type and their equilibrium states by 
using the Laplace transform (Hirata~\cite{Hir1}, Coelho~\cite{Coe} and Collet) or using the moment method
(Pitskel~\cite{Pit} and others as for instance in~\cite{H2} for rational maps). 
Galves and Schmitt~\cite{GS} showed that entry times are exponentially 
distributed for $\phi$-mixing systems and also gave error terms by a method that
later was considerable expanded and sharpened by Abadi~\cite{Aba} for $\phi$-mixing 
and later even some $\alpha$-mixing systems.
A more elementary counting approach by Abadi and Vergne~\cite{AV3} proved that 
$\phi$-mixing measures have exponentially decaying entry and return times.
They also have results on multiple return times which are Poissonian for $\psi$-mixing measures.
In~\cite{BSTV} a reduction to induced maps was used to obtain the limiting
distribution of return times. The important connection there was a result 
that establishes that the limiting return times distribution for the given system is
the same as the limiting return times distribution of an induced system. 
This is the result we expand and sharpen in this note.
It allows to obtain distribution results for many more systems because instead
of checking mixing properties for the given system one can consider the 
jump transform which is an induced map on a suitable subset so that
 the induced map has good expanding properties that result in 
 exponential decay of correlations. Such an approach has been particularly 
 successful in the study of parabolic interval maps like the Manneville-Pommeau
 map which is non-uniformly expanding because it has a parabolic fixed point
 where the derivative is equal to $1$.  The induced map on any interval not including
 the parabolic fixed point is uniformly expanding and its statistical properties
 can more easily be  analysed exploiting the quasi compactness of the transfer 
 operator and the consequential exponential decay of correlations.

\section{Return times and the induced map}

Let $\Omega$ be a measure space with $\sigma$-algebra $\mathscr{B}$ and 
$\mu$ a probability measure. Moreover $T$ is a measure preserving map on $\Omega$.
We assume that $\mu$ is ergodic.
For $U\subset\Omega$ we define the function
$$
\tau_U(x)=\min\{j\ge1: T^jx\in U\},
$$
which is the {\em entry time} for $x$ if $x\in\Omega$ and is the {\em return time} if $x\in U$.
We put $\tau_U(x)=\infty$ if the forward orbit of $x$ never enters $U$.
 Poincar\'e's recurrence theorem states that $\tau_U(x)<\infty$ for almost every $x\in U$ 
 for any finite $T$-invariant measure $\mu$ on $\Omega$ and Kac's theorem~\cite{Kac} tells us
 that $\tau_U$ is integrable on $U$. In fact 
$$
\int_U\tau_U(x)\,d\mu(x)=1
$$
if $\mu(U)>0$. Let us note that $\tau_U$ is not necessarily integrable over $\Omega$,
in fact $\int_\Omega\tau_U\,d\mu<\infty$ if and only if $\tau_U$ is square integrable over $U$.

For a subset $U\subset\Omega$, $\mu(U)>0$, let us denote by 
$\hat{T}=T^{\tau_U}:U\rightarrow U$ the {\em induced map}. $\hat{T}$ exists
by Poincar\'e's (or Kac's) theorem almost everywhere. We also have
the {\em induced measure} $\hat\mu$ which is defined on $U$ by
$\hat\mu(A)=\frac{\mu(A)}{\mu(U)}$ for all measurable $A\subset U$.
The induced measure $\hat\mu$ is $\hat{T}$-invariant and ergodic (see e.g.~\cite{Pet}).

\subsection{Entry times distributions}
Let $B\subset \Omega$ ($\mu(B)>0$) and put for (parameter values) $t>0$
$$
F_B(t)=\mathbb{P}\left(\tau_B>\frac{t}{\mu(B)}\right)
=\mu\left(\left\{x\in\Omega: \tau_B(x)>\frac{t}{\mu(B)}\right\}\right)
$$
for the entry time distribution to $B$.
The entry times distribution $F_B(t)$ is locally constant on intervals of length $\mu(B)$
 and has jump discontinuities at values $t$ which are integer multiples of $\mu(B)$. 
 For any $s\in\mathbb{N}_0$ one has
\begin{equation}\label{regularity}
\{\tau_B>s+1\}=T^{-1}\{\tau_B>s\}\setminus T^{-1}B
\end{equation}
and consequently
 $$
 \mathbb{P}(\tau_B=s+1)=\mathbb{P}(\tau_B>s)-\mathbb{P}(\tau_B>s+1)\le\mu(B)
 $$
which shows that the jumps at the discontinuities are at most $\mu(B)$. Hence
$$
|F_B(t)-F_B(s)|\le|t-s|+\mu(B)
$$
for all $t,s>0$.

Now let $B_j\subset\Omega$ ($\mu(B_j)>0$) be a sequence of subsets so that 
$\mu(B_j)\rightarrow0^+$ as $j\rightarrow\infty$. We want to assume that 
$F_{B_j}(t)$ converges pointwise (in $t$) to a limiting distribution $F(t)$ as $j\rightarrow\infty$. 
 Let us note that the regularity~(\ref{regularity}) the limiting distribution $F(t)$ is Lipschitz continuous 
 with Lipschitz constant $1$ and consequently $F$ is  continuous.

Lacroix~\cite{Lac} has shown that if $F(t)$ is an eligible limiting distribution, that is 
it satisfies $F(0)=1$, is continuous, convex, monotonically decreasing on $(0,\infty)$ and 
$F(t)\rightarrow0^+$ as $t\rightarrow\infty$,
then for any ergodic $T$-invariant probability measure $\mu$ there exists a sequence of 
positive measure sets $B_j\subset\Omega$ so that $\mu(B_j)\rightarrow0$ and
such that $F(t)=\lim_{j\rightarrow\infty}F_{B_j}(t)$ for every $t\in(0,\infty)$.
The sets $B_j$ are typically pretty wild looking and in particular won't be topological
balls or cylinder sets for a given partition. A similar result was shown for return times in~\cite{KL}
although the two results are equivalent by~\cite{HLV}.

  For a positive measure subset $U\subset\Omega$ let us now consider the induced system
  $(U,\hat{T},\hat\mu)$ which carries the entry time function $\hat\tau_B(x)>\min\{j\ge1:\hat{T}^jx\in B\}$
  for sets $B\subset U$, $\hat\mu(B)>0$. As above we can then define the entry times distribution
  $$
  \hat{F}_{B}(t)=\hat{\mathbb{P}}\left(\hat\tau_{B}>\frac{t}{\hat\mu(B)}\right)
  =\hat\mu\left(\left\{x\in U: \hat\tau_B(x)>\frac{t}{\hat\mu(B)}\right\}\right)
  $$

\vspace{3mm}

\noindent
The following theorem shows that a restricted system $(U,\hat{T},\hat\mu)$ has the same
limiting entry times distribution as the original system $(\Omega, T,\mu)$.

\begin{theorem}
 Let $\mu$ be ergodic, $U\subset\Omega$, $\mu(U)>0$. 
 Assume there exists a sequence of sets $B_j\subset U$, $\mu(B_j)\rightarrow0^+$, so that
 either the limiting entry times distribution for $(\Omega,T,\mu)$
$$
 F(t)=\lim_{j\rightarrow\infty}F_{B_j}(t)
 $$
 exists, or the limiting entry times limiting distribution
 $$
 \hat{F}(t)=\lim_{j\rightarrow\infty}\hat{F}_{B_j}(t)
 $$
 for the induced system exists $(U,\hat{T},\hat\mu)$ exists.
 
 Then both limiting entry times distributions exsit and moreover $F(t)=\hat{F}(t)$ for all $t>0$.
\end{theorem}

  \noindent {\it Proof.} Let $B=B_j$. We first relate $\tau_B$ to $\hat\tau_B$ ($B\subset U,\mu(B)>0$).
 If we put $m=\hat\tau_B(x)$, $x\in U$, then
 $$
 \tau_B(x)=\tau_U(x)+\tau_U(\hat{T}x)+\tau_U(\hat{T}^2x)+\cdots+\tau_U(\hat{T}^{m-1}x)=n^m(x),
 $$
 where we wrote the ergodic sum of the function $n=\tau_U|_U$ for the return time on $(U,\hat{T})$.
 By the Birkhoff ergodic theorem on $(U,\hat{T},\hat\mu)$ we get as $\hat\mu$ is ergodic:
 $$
 \frac1m\tau_B(x)=\frac1mn^m(x)\rightarrow\int_Un(x)\,d\mu(x)
 =\int_U\tau_U(x)\,\frac{d\mu(x)}{\mu(U)}=\frac1{\mu(U)}
 $$
 as $m\rightarrow\infty$ by Kac's theorem for almost every $x\in U$.
 
 Let $\varepsilon>0$, then there exists $G_\varepsilon\subset U$, and $M_\varepsilon\in\mathbb{N}$
 so that 
 $$
 \left|\frac1mn^m(x)-\frac1{\mu(U)}\right|<\varepsilon
 \hspace{10mm} \forall\;x\in G_\varepsilon,\;m\ge M_\varepsilon
 $$
 and $\mu(G_\varepsilon^c)<\varepsilon$. Thus 
 $$
 \tau_B(x)=\sum_{j=0}^{\hat\tau_B(x)-1}\tau_U\circ \hat{T}^j=\frac{\hat\tau_B(x)}{\mu(U)}+\mathcal{O}(\hat\tau_B(x)\varepsilon)
 $$
 for all $x\in G_\varepsilon$ such that $\hat\tau_B(x)\ge M_\varepsilon$. 
 Since $\tau_U$ is integrable on $U$ there exists a $\delta>0$ (depending on $\varepsilon$)
 so that $\int_S\tau_U\,d\mu<\varepsilon$ for any set $S\subset U$ for which $\mu(S)<\delta$.
 We can assume that $\mu(G_\varepsilon^c)<\min(\delta,\varepsilon)$.

 For $j=0,1,2,\dots$ put $ A_j=\Omega\setminus T^{-j}U=T^{-j}U^c$ and 
  $$
  D_j^k=\bigcap_{\ell=j}^k A_\ell=\{x\in\Omega: T^\ell x\not\in U\;\forall \ell=j,\dots,k\}
  $$
  for $0\le j\le k$. Then for any $j\in\mathbb{N}$
  $$
  \{x\in\Omega: \tau_U(x)=j\}=T^{-j}U\cap D_1^{j-1}=D_1^{j-1}\setminus D_1^j.
  $$
  On the other hand we also have 
  $$
\left\{x\in U: \tau_U(x)\ge j\right\}=U\cap D_1^{j-1}.
  $$
 We now do the following decomposition (as $D_1^{j-1}=\{x\in\Omega: \tau_U(x)\ge j\}$):
 \begin{eqnarray*}
 F_B(t) &=&\int_\Omega\chi_{\tau_B>s}\,d\mu\\
 &=&\sum_j\int_{\{\tau_U=j\}}\chi_{\tau_B>s}\,d\mu\\
  &=&\sum_j\left(\int_{D_1^{j-1}}\chi_{\tau_B>s}\,d\mu
  -\int_{D_1^{j}}\chi_{\tau_B>s}\,d\mu\right),
 \end{eqnarray*}
as $D_1^{j}\subset D_1^{j-1}$, where we wrote $s=\frac{t}{\mu(B)}$.
For the second term in the last line consider
 $$
 \int_{D_0^{j-1}}\chi_{\tau_B>s}\,d\mu
 =\int_\Omega \left(\chi_{D_0^{j-1}}\chi_{\tau_B>s}\right)\circ T\,d\mu
 = \int_{D_1^{j}}\chi_{\tau_B>s}\circ T\,d\mu
 $$
 as $T^{-1}D_0^{j-1}=D_1^j$. The inclusions
 $$
 \{\tau_B> s+1\}\subset T^{-1}\{\tau_B> s\}\subset\{\tau_B> s+1\}\cup T^{-1}B
 $$
imply the inequalities
 $$
 \int_{D_1^{j}}\chi_{\tau_B> s+1}\,d\mu
 \le\int_{D_1^{j}}\chi_{\tau_B>s}\circ T\,d\mu
 \le\int_{D_1^{j}}\chi_{\tau_B>s+1}\,d\mu
 +\int_{D_1^{j}}\chi_{T^{-1}B}\,d\mu
 $$
 where the last integral is equal to zero as $T^{-1}B\cap D_1^j=T^{-1}(B\cap D_0^{j-1})=\emptyset$
 because $D_0^{j-1}\subset U^c\subset B^c$ for $j\ge1$. Thus
 $$
 \int_{D_1^{j}}\chi_{\tau_B>s}\,d\mu
 =  \int_{D_1^{j}}\chi_{\tau_B>s-1}\circ T\,d\mu
 = \int_{D_0^{j-1}}\chi_{\tau_B>s-1}\,d\mu
 $$
 and 
\begin{eqnarray*}
 F_B(t)&=&\sum_j\left( \int_{D_1^{j-1}}\chi_{\tau_B> s}\,d\mu- \int_{D_0^{j-1}}\chi_{\tau_B> s-1}\,d\mu\right)\\
&=& \sum_j\int_{C_j}\chi_{\tau_B> s}\,d\mu+ E_0\\
&=&\int_U\tau_U\chi_{\tau_B> s}\,d\mu+ E_0
\end{eqnarray*}
 where we put $C_j=D_1^{j-1}\setminus D_0^{j-1}=\{x\in U: \tau_U(x)\ge j\}$ and used
 that $\sum_j\chi_{C_j}=\tau_U$ on $U$ (as $|\{j: x\in C_j\}|=\tau_U(x)\forall x\in U$).
 To estimate the error term $E_0$ note that 
 $\mathbb{P}(\tau_B=s)=\mathbb{P}(\tau_B>s-1)-\mathbb{P}(\tau_B>s)\le\mu(B)$
 (see the remark preceding the theorem).
 Since $D_0^{j-1}=\{x\in U^c:\tau_U(x)\ge j\}$ one has
$$
-E_0=\sum_j\int_{D_0^{j-1}}\chi_{\tau_B=s}\,d\mu
\le\int_\Omega\tau_U\chi_{\tau_B=s}\,d\mu.
$$
In order to show that $E_0$ goes to zero as $\mu(B)$ decreases to zero let us put 
$B_{k,j}=U_k\cap T^{-j}\{\tau_B=s\}$ for $j=0,1,\dots,k-1$ and $k=1,2,\dots$, 
where $U_k=\{x\in U: \tau_U(x)=k\}$.
Then $B_{k,j}\cap B_{k,i}=\emptyset$ if $i\not=j$ because if there were an $x\in B_{k,j}\cap B_{k,i}$
it would imply that $T^jx, T^ix\in\{\tau_B=s\}$ and therefore, assuming $j>i$, one would get
the contradiction $s=\tau_B(T^jx)>\tau_B(T^ix)=s$.
Hence the sets $B_{k,j}$ are 
pairwise disjoint for $k\in\mathbb{N}$ and $j=0,\dots,k-1$. 
In other words, for every $x\in B_k$ there is a unique $j\in[0,k)$ so that $T^jx\in\{\tau_B=s\}$.
Consequently $\{\tau_B=s\}=\dot{\bigcup}_{k=1}^\infty\dot{\bigcup}_{j=0}^{k-1}T^jB_{k,j}$,
and since $\mu(T^jB_{k,j})\ge\mu(B_{k,j})$ we obtain
$$
\mu(\tilde{B})=\sum_{k=1}^\infty\sum_{j=0}^{k-1}\mu(B_{k,j})
\le\sum_{k=1}^\infty\sum_{j=0}^{k-1}\mu(T^jB_{k,j})
=\mathbb{P}(\tau_B=s)\le\mu(T^{-s}B)=\mu(B),
$$
where $\tilde{B}=\dot{\bigcup}_{k=1}^\infty\dot{\bigcup}_{j=0}^{k-1}B_{k,j}$ ($\tilde{B}\subset U$).
For $x\in \tilde{B}$ for which $T^jx\in\{\tau_B=s\}$  one has 
$\tau_U(x)\ge\tau_U(T^jx)$ and therefore
$$
|E_0|\le \int_{\tilde{B}}\tau_U\,d\mu<\varepsilon
$$
 as we can assume that $\mu(B)<\delta$.

Replacing $\tau_B$ by $\hat\tau_B$ using the relation 
$\frac{\tau_B}{\hat\tau_B}=\frac1{\mu(U)}+\mathcal{O}(\varepsilon)$ implies
$\tau_B=\frac{\hat\tau_B}{\mu(U)}+\eta$ where $\eta:U\rightarrow\mathbb{R}$ has the bound
$|\eta|\le\varepsilon\mu(U)\hat\tau_B\le \varepsilon\hat\tau_B$. Hence
 $$
 F_B(t)=\int_U\tau_U \chi_{\hat\tau_B>\frac{t}{\hat\mu(B)}+\eta}\,d\mu +E_0.
 $$ 
 Now we want to introduce a power $k$ of the induced map $\hat{T}$ so that we can average over $k$
 and use the ergodic theorem on $(U,\hat{T}, \hat\mu)$. By $\hat{T}$-invariance of $\hat\mu$
\begin{eqnarray*}
 F_B(t)&=&\int_U\left(\tau_U
 \chi_{\hat\tau_B>\frac{t}{\hat\mu(B)}+\eta}\right)\circ\hat{T}^k\,d\mu+E_0\\
 &=&\int_{G_\varepsilon}\left(\tau_U
 \chi_{\hat\tau_B>\frac{t}{\hat\mu(B)}+\eta}\right)\circ\hat{T}^k\,d\mu+E_0+H_k,
\end{eqnarray*}
 where we get for the error
 $$
  H_k=\int_{G_\varepsilon^c}\left(\tau_U \chi_{\tau_B>\frac{t}{\mu(B)}}\right)\circ \hat{T}^k\,d\mu
  \le\int_{G_\varepsilon^c}\tau_U\circ\hat{T}^{k}\,d\mu
  =\int_{\hat{T}^{-k}G_\varepsilon^c}\tau_U\,d\mu<\varepsilon
 $$
since by assumption $\mu(\hat{T}^{-k}G_\varepsilon^c)=\mu(G_\varepsilon^c)<\delta$ as $\mu$
restricted to $U$ is $\hat{T}$-invariant.

 In the principal term we want to exploit the identity $\sum_j\chi_{C_j}=\tau_U$ 
on $U$. For that purpose note that
 $$
 \left\{x\in U:\hat\tau_B(\hat{T}^{k}x)\ge s\right\}\setminus\bigcup_{\ell=1}^{k-1}\hat{T}^{-\ell}B
 =\left\{x\in U:\hat\tau_B(x)\ge s+k\right\}
 $$
 which yields (here we use $s=\frac{t}{\hat\mu(B)}+\eta$)
 $$
 F_B(t)=\int_{G_\varepsilon}\left(\tau_U
  \chi_{\hat\tau_B>\frac{t}{\hat\mu(B)}+\eta+k}\right)\circ\hat{T}^{k}\,d\mu
 + E_0+H_k+K_k,
 $$
 where the individual errors are bounded by:
 $$
 K_k\le\int_{G_\varepsilon}\tau_U\circ\hat{T}^k\sum_{\ell=1}^{k-1}\chi_B\circ\hat{T}^\ell\,d\mu.
 $$
 We now estimate the average error over $k\in\{0,1,\dots,n-1\}$:
 \begin{eqnarray*}
 \hat{K}_n&=&\frac1n\sum_{k=0}^{n-1}K_k\\
 &=&\int_{G_\varepsilon}\frac1n\sum_{k=0}^{n-1}\sum_{\ell=1}^{k-1}
 (\tau_U\circ\hat{T}^k)(\chi_B\circ\hat{T}^\ell)\,d\mu\\
  &=&\int_{G_\varepsilon}\sum_{\ell=1}^{n-2}(\chi_B\circ\hat{T}^\ell)
\left(\frac1n\sum_{k=\ell+1}^{n-1}\tau_U\circ\hat{T}^k\right)\,d\mu\\
&\le&c_1\frac1{\mu(U)}\int_{G_\varepsilon}\sum_{\ell=1}^{n-2}(\chi_B\circ\hat{T}^\ell)\,d\mu\\
&\le&c_1n\hat\mu(B)
 \end{eqnarray*}
 where we used the estimate
 $$
 \frac1n\sum_{k=\ell+1}^{n-1}\tau_U\circ\hat{T}^k
 \le\frac1n\sum_{k=0}^{n-1}\tau_U\circ\hat{T}^k\le\frac1{\mu(U)}+\varepsilon
 \le c_1\frac1{\mu(U)}
 $$
 for some constant $c_1$ and for all $x\in G_\varepsilon$ provided $n\ge M_\varepsilon$.
Thus
 $$
 F_B(t)=\frac1n\sum_{k=1}^n\int_{G_\varepsilon}\left(\tau_U\circ\hat{T}^k\right)
 \chi_{\hat\tau_B>\frac1{1-\eta'}\frac{t+k\hat\mu(B)}{\hat\mu(B)}}\,d\mu
 + E_0+\hat{H}_n+\hat{K}_n,
 $$
 where $\hat{H}_n=\frac1n\sum_{k=0}^{n-1}H_k<\varepsilon$
 and $\eta'':U\rightarrow\mathbb{R}$ satisfies the bound $|\eta''|<\varepsilon$.
Consequently (as $|E_0+\hat{K}_n|<\varepsilon+c_1n\hat\mu(B)$)
 \begin{eqnarray*}
F_B(t)&=&\int_{G_\varepsilon}\frac1n\sum_{k=1}^n\left(\tau_U\circ\hat{T}^k\right)
 \chi_{\hat\tau_B>(1+\eta'')\frac{t}{\hat\mu(B)}}\,d\mu
+\mathcal{O}(\varepsilon+n\hat\mu(B))\\
&=&\int_{G_\varepsilon}
 \chi_{\hat\tau_B>(1+\eta'')\frac{t}{\hat\mu(B)}}\,d\hat\mu
+\mathcal{O}(\varepsilon+n\hat\mu(B))
 \end{eqnarray*}
 as $\frac1n\sum_{k=1}^n\tau_U\circ\hat{T}^k=\frac1{\mu(U)}+\mathcal{O}(\varepsilon)$ on $G_\varepsilon$,
 where $\eta'':U\rightarrow\mathbb{R}$ satisfies $|\eta''|<c_2|\eta'|+\frac{n}t\hat\mu(B)$ ($c_2>0$).
 To adjust for the `time shift' in the lower bound of the entry function, we use the fact that 
 $|\hat{F}_B(t)-\hat{F}_B(s)|\le |t-s|+\hat\mu(B)$ and thus obtain (for a $c_3$)
 $$
 |F_B(t)-\hat{F}_B(t)|<c_3\varepsilon+\left(\frac{n}t+1\right)\hat\mu(B)+c_1n\hat\mu(B)
 $$
 for all $n\ge M_\varepsilon$. If $\mu(B_j)$ is small enough so that
  $\hat\mu(B_j)<\min(\frac{\varepsilon t}{M_\varepsilon},\frac\varepsilon{(c_1+1)n})$
  (if we choose $n=M_\varepsilon$ this requires 
  $\hat\mu(B_j)<\frac\varepsilon{M_\varepsilon}\min(t,\frac1{c_1+1})$)
 then
  $$
 |F_{B_j}(t)-\hat{F}_{B_j}(t)|<(c_3+1)\varepsilon
 $$
 and as $\mu(B_j)\rightarrow0$ ($j\rightarrow\infty$) we obtain 
 $|F(t)-\hat{F}(t)|<(c_3+1)\varepsilon$ for any positive $\varepsilon$.
Thus if the limiting distribution $F(t)$ exists then also the limiting distribution $\hat{F}(t)$ exists and 
vice versa. Moreover we obtain equality: $F(t)=\hat{F}(t)$ for all $t>0$.
\qed 
 
\subsection{Return times distributions}
  The restriction of the function $\tau_B$ to the set $B\subset\Omega$ is 
 called the {\em return time function} and  we correspondingly call
 $$
 \tilde{F}_B(t)=\mathbb{P}_B\left(\tau_B>\frac{t}{\mu(B)}\right)
 $$
 the {\em return times distribution}.
 For instance, if $\Omega$ is the shiftspace $\Sigma$
 and $B=U(x_0x_1\cdots x_{n-1})$ is an $n$-cylinder then $\tau_B(\vec{x})$ for $\vec{x}\in B$
 measures the `time' it takes to see the word $x_0x_1\cdots x_{n-1}$ again, that is
 $$
 \tau_B(\vec{x})=\min\{j\ge1: x_jx_{j+1}\cdots x_{j+n-1}=x_0x_1\cdots x_{n-1}\}.
 $$
 The function $\tilde{F}_B(t)$ then measures the probability to see the first $n$-word
 again after rescaled time $t/\mu(B)$.
 
 Similarly for the induced system $(U,\hat{T},\hat\mu)$ we have the return times distribution
 $$
 \hat{\tilde{F}}_B(t)=\hat{\mathbb{P}}_B\left(\hat\tau_B>\frac{t}{\hat\mu(B)}\right)
 =\hat\mu\left(\left\{x\in B: \hat\tau_B(x)>\frac{t}{\hat\mu(B)}\right\}\right).
 $$
 In order to get a similar result on the relation between return times for the original system
 and the induced system, we will need the following result.
  
 \begin{proposition} \cite{HLV}
 Let $B_j\subset\Omega$ ($\mu(B_j)>0$) be a sequence of sets so that $\mu(B_j)\rightarrow0^+$.
 If one of the limits $F(t)=\lim_{j\rightarrow\infty} F_{B_j}(t)$, $ \tilde{F}(t)=\lim_{j\rightarrow\infty}  \tilde{F}_{B_j}(t)$
 exists (pointwise) then so does the other limit and moreover
 $$
F(t)=\int_t^\infty  \tilde{F}(s)\,ds.
 $$
 \end{proposition}
 
 \noindent While the limiting entry times distribution $F(t)$ is always Lipshitz continuous with
 Lipshitz constant $1$, the same does not apply to the limiting return times distribution 
 $\tilde{F}(t)$ which in fact can have (at most countable many) discontinuities.
 In particular, if the sets $B_j$ contract to a periodic point, then $\tilde{F}(t)$ will have a
 discontinuity at $t=0$ with $\lim_{t\rightarrow0^+}\tilde{F}(t)<1$.
Also note that since the limiting entry distribution $F$ is Lipschitz continuous
  the limiting return distribution $\tilde{F}(t)$ is monotonically decreasing to zero
which implies that $F(t)$ is in fact always convex.
 
 One consequence of this result is that the limiting entry times distribution and return times
distribution are the same only if they are exponential, that is $\tilde{F}=F$ if only if
$F(t)=\tilde{F}(t)=e^{-t}$. We use this proposition to obtain the corresponding result of Theorem~1
for the limiting return times distribution.

\begin{theorem}
 Let $\mu$ be ergodic, $U\subset\Omega$, $\mu(U)>0$. 
 Assume there exists a sequence of sets $B_j\subset U$, $\mu(B_j)\rightarrow0^+$, so that 
one of the two limiting return times distribution
$$
\mbox{either}\hspace{5mm} \tilde{F}(t)=\lim_{j\rightarrow\infty} \tilde{F}_{B_j}(t),\hspace{5mm}
\mbox{or}\hspace{5mm} \hat{ \tilde{F}}(t)=\lim_{j\rightarrow\infty}\hat{ \tilde{F}}_{B_j}(t)
 $$ 
 exists.
 
 Then both limiting return times distributions exist and 
 moreover at every point of continuity $t\in\mathbb{R}^+$ one has equality
$\tilde{F}(t)=\hat{ \tilde{F}}(t)$.
\end{theorem}

\noindent This is the result that was proven in~\cite{BSTV} in 2003 for Radon measures on Riemann 
manifolds using the Lebesgue Density theorem. The limit there was along metric balls 
$B_j$ that shrink to a point $x$  and with the implication that the 
existence of the limiting return times distribution in the induced system (plus the non-degeneracy condition
 $\tilde{F}(0^+)=1$) implies the limiting return times distribution for the entire system
 and that the two limiting distributions are equal at points of continuity.
 
 \vspace{3mm}
 
\noindent {\it Proof of Theorem~3.} Assume that, say, the limit 
 $ \tilde{F}(t)=\lim_{j\rightarrow\infty} \tilde{F}_{B_j}(t)$ exists. By Proposition~2 this implies
 the also the limiting distribution $F (t)=\lim_{j\rightarrow\infty} F_{B_j}(t)$ exists. 
 By Theorem~1 we get that the limit $\hat{F}(t)=\lim_{j\rightarrow\infty}  \hat{F}_{B_j}(t)$
 exists and satisfies $\hat{F}=F$. Again by Proposition~2 this implies the limit
  $\hat{\tilde{F}}(t)=\lim_{j\rightarrow\infty} \hat{\tilde{F}}_{B_j}(t)$ exists. Thus,
  since
  $$
  \int_t^\infty\tilde{F}(s)\,d\mu(s)=F(t)=\hat{F}(t)=  \int_t^\infty\hat{\tilde{F}}(s)\,d\mu(s)
  $$
  for all $t>0$ we conclude that $\tilde{F}(t)=\hat{\tilde{F}}(t)$ at all points $t$ of
  continuity.
  
  Similarly on shows that the limit  $\hat{\tilde{F}}(t)=\lim_{j\rightarrow\infty} \hat{\tilde{F}}_{B_j}(t)$
  implies the return times limiting distribution $ \tilde{F}(t)=\lim_{j\rightarrow\infty} \tilde{F}_{B_j}(t)$
  for the whole system and also equality of the limiting distributions $\tilde{F}(t)=\hat{\tilde{F}}(t)$
   at points of continuity.
  \qed
  
  \vspace{3mm}
  
  \noindent Let us note that since sole the requirements in Theorem~1 and~2 are the existence 
  of the limits they apply in particular also to examples of Lacroix and Kupsa~\cite{Lac,KL} 
  where for any ergodic transformation they produce a sequence $B_j$ that realises an
  arbitrary given (eligible) limiting distribution.
  
\subsection{Example} Here we give an example where it is easy to find the limiting 
entry/return times distributions for the induced map.
   We consider the shift space $\Omega=\mathbb{N}^\mathbb{Z}$ with the shift transformation $\sigma$.
 To define the invariant measure $\mu$ we give on the state space $\mathbb{N}$ 
  the transition probabilities: Let $p_i\in(0,1), i=1,2,\dots$, be a sequence, then
  we allow for the transition $i\rightarrow i+1$ with probability $p_i$ and 
  for the transition $i\rightarrow 1$ with probability $q_i=1-p_i$. In other words,
  we can define a stochastic matrix $M$ by
  $$
  \left\{\begin{array}{rcl} 
  M_{j,1}&=&q_j\\
  M_{j,j+1}&=&p_j\\
  M_{j,k}&=&0\mbox{ otherwise, i.e.\ if $k\not=1$ or $k\not=j+1$}
  \end{array}\right.,
  $$
  where the transition probability of the transition $j\rightarrow k$ is given by the entry $M_{j,k}$.
  Then
  $M\mathds{1}=\mathds{1}$ as $\sum_{k=1}^\infty M_{j,k}=M_{j,1}+M_{j,j+1}=q_j+p_j=1\forall j$
  and $M$ has the left eigenvector $\vec{x}=(x_1,x_2,\dots)$ (for the dominant eigenvalue $1$)
  which satisfies
  \begin{eqnarray*}
  q_1x_1+q_2x_2+q_3x_3+\cdots&=&x_1\\
  x_jp_j&=&x_{j+1}\;\;\mbox{ for $j=1,2,\dots$}.
  \end{eqnarray*}
One sees that the components of the left eigenvector are  $x_j=x_1P_j, j=2,3,\dots$,
  where $P_j=\prod_{i=1}^{j-1}p_i$ ($P_1=1$) and $x_1$ is chosen to  make $\vec{x}$ a 
  probability vector ($x_1^{-1}=\sum_jP_j$). We assume $x_1>0$. The first equation above is
  satisfied as $\sum_jq_jx_j=x_1\sum_j(P_j-P_{j+1})=x_1$ if $P_j\rightarrow0$ as
  $j\rightarrow\infty$.
  In this way we obtain a shift invariant probability measure $\mu$ on $\Omega$
  which is ergodic as one can go from any state $i$ to any other state $j$ with positive
  probability.
  
  Put $A_j=\{\vec\omega\in\Omega:\omega_0=j\}, j=1,2,\dots$, and let 
  $U=A_1$ be the return set with return/entry time function $\tau_U$.
If we put $A_{j,k}=A_j\cap \{\tau_U=k\}$ then $\vec\omega\in A_{j,k}$ is of the form
$\omega_0\omega_1\cdots\omega_k=j(j+1)(j+2)\cdots(j+k-2)(j+k-1)1$
(symbol sequence of length $k+1$).
One has 
$$
\mu(A_{j,k})=\mu(A_j)p_jp_{j+1}\cdots p_{j+k-2}q_{j+k-1}
=x_1P_j\frac{P_{j+k-1}}{P_j}q_{j+k-1}
=x_1P_{j+k-1}q_{j+k-1}
$$
as $\mu(A_j)=x_j=x_1P_j$. Let $\mathcal{D}$ be the countably infinite partition of $U$ 
whose partition elements are $D_j=\{\omega\in U:\tau_U(\omega)=j\}$ ($D_j=A_{1,j}$).
 The induced map 
$\hat\sigma:U\to U$ is a Bernoulli shift on $\hat\Omega=\mathcal{D}^\mathbb{Z}$
and the induced measure $\hat\mu$ is the Bernoulli measure with weights  
$\hat\mu(D_j)=\frac1{\mu(U)}x_1q_jP_j$, where $\mu(U)=\sum_jx_1q_jP_j$. If we denote by $B_n$ the $n$-cylinder which contains
a given point $\hat\omega\in\hat\Omega$ then the entry times $\hat F_{B_n}(t)$ converge
to the exponential distribution $e^{-t}$ as $n\rightarrow\infty$ for almost every $\hat\omega$.
Hence we conclude that entry times $F_{B_n}$ for the map $\sigma$ on $\Omega$
also converge to the limiting distribution $e^{-t}$ almost surely.

\vspace{3mm}

 \noindent {\bf Remark.} Kac's  theorem states that  the return time function $\tau_U$
  is integrable over  $U$ and  also gives the value  of the integral. We can use this
  example to achieve that $\tau_U$ is not integrable over the entire space $\Omega$
 alhough the measure is ergodic. The integral of $\tau_U$ over the entire space is
$$
\int_\Omega\tau_U\,d\mu=\sum_{j,k}k\mu(A_{j,k})=\sum_{j,k}kx_1P_{j+k-1}q_{j+k-1}.
$$
If we choose $p_i=\left(\frac{i}{i+1}\right)^\alpha$ for some $\alpha\in(1,2)$ then
$P_j=\prod_{i=1}^{j-1}\left(\frac{i}{i+1}\right)^\alpha=\frac1{j^\alpha}$ and since
the $P_j$ are summable, $x_1=\left(\sum_jP_j\right)^{-1}$ is well defined and
positive. Then
\begin{eqnarray*}
\int_\Omega\tau_U\,d\mu&=&x_1\sum_kk\sum_j\frac1{(j+k-1)^\alpha}q_{j+k-1}\\
&\ge&c_1x_1\sum_kk\sum_j\frac1{(j+k-1)^{\alpha+1}}\\
&\ge&c_2\sum_k\frac{k}{k^\alpha}=\infty,
\end{eqnarray*}
as $\alpha<2$, where we used that $q_{j+k-1}=1-\left(1-\frac1{j+k-1}\right)^\alpha\ge c_1 \frac1{j+k-1}$
for some $c_1>0$. We thus see that the integral of $\tau_U$ over the entire space $\Omega$ diverges.

This can be converted to  an example on a two-state shiftspace $\Sigma\subset\{0,1\}^\mathbb{Z}$
 by the single element mapping $\pi:\Omega\rightarrow\Sigma$ which maps $\pi(1)=1$
 and collapses all other symbols to $0$, i.e. $\pi(j)=0, j=2,3,\dots$. The measure $\mu$
 is sent to the probability measure $\nu=\pi^*\mu$ which is invariant under the shift map.
 
 In fact  $\int_\Omega\tau_U\,d\mu$ is finite if and only if $\int_U\tau_U^2\,d\mu$ is 
 finite. So the above example is an example where the return time to $U$ is not square integrable
 over $U$.



\end{document}